\documentclass[intlimits,12pt,reqno]{amsart}
\usepackage{latexsym,amsthm,amsmath,amssymb,mathrsfs,mathtools}
\usepackage[T1]{fontenc}
\usepackage[utf8]{inputenc}
\usepackage[mathscr]{eucal}
\usepackage{enumerate}
\usepackage{enumitem}

\let\oldmarginpar\marginpar
\renewcommand{\marginpar}[2][rectangle,draw,text width= 2cm,rounded corners]{
    \oldmarginpar{
    \scriptsize \tikz \node at (0,0) [#1]{#2};}
    }
\usepackage{tikz}
\usetikzlibrary{arrows,calc,patterns,cd}
\usepackage{wrapfig}

\usepackage{hyperref}
\usepackage{xcolor}
\hypersetup{
    colorlinks=true,
    linkcolor=blue,
    citecolor=blue,
    urlcolor=blue
}

plus 6pt minus 12pt
\def\mvint_#1{\mathchoice
          {\mathop{\vrule width 6pt height 3 pt depth -2.5pt
                  \kern -9pt \intop}\limits_{\kern -3pt #1}}%
          {\mathop{\vrule width 5pt height 3 pt depth -2.6pt
                  \kern -6pt \intop}\nolimits_{#1}}%
          {\mathop{\vrule width 5pt height 3 pt depth -2.6pt
                  \kern -6pt \intop}\nolimits_{#1}}%
          {\mathop{\vrule width 5pt height 3 pt depth -2.6pt
                  \kern -6pt \intop}\nolimits_{#1}}}

\def\vi{\varphi}

\newcommand{\MM}{\mathcal M}

\newcommand{\bbbr}{\mathbb R}

\newcommand{\bbbn}{\mathbb N}

\newcommand{\N}{\mathbb N}

\newcommand{\R}{\mathbb R}
\newcommand{\overbar}[1]{\mkern 1.7mu\overline{\mkern-1.7mu#1\mkern-1.5mu}\mkern 1.5mu}

\newcommand{\eps}{\varepsilon}

\newtheorem{theorem}{Theorem}[section]
\newtheorem*{theorem*}{Theorem}
\newtheorem{lemma}[theorem]{Lemma}
\newtheorem{corollary}[theorem]{Corollary}
\newtheorem{proposition}[theorem]{Proposition}

\theoremstyle{definition}

\newtheorem{remark}[theorem]{Remark}
\newtheorem*{remark*}{Remark}

\newcommand{\supp}{{\rm supp}\,}

\DeclareMathOperator{\Lip}{Lip}
\DeclareMathOperator{\inte}{int}
\DeclareMathOperator{\dist}{dist}

\newcommand*{\loc}{{\mathrm{loc}}}

\usepackage{setspace}

\makeatletter
\renewcommand{\tocsection}[3]{%
  \indentlabel{\@ifnotempty{#2}{\bfseries\ignorespaces#1 #2\quad}}\bfseries#3}
\renewcommand{\tocsubsection}[3]{%
  \indentlabel{\@ifnotempty{#2}{\ignorespaces#1 #2\quad}}#3}

\newcommand\@dotsep{4.5}
\def\@tocline#1#2#3#4#5#6#7{\relax
  \ifnum #1>\c@tocdepth 
  \else
    \par \addpenalty\@secpenalty\addvspace{#2}%
    \begingroup \hyphenpenalty\@M
    \@ifempty{#4}{%
      \@tempdima\csname r@tocindent\number#1\endcsname\relax
    }{%
      \@tempdima#4\relax
    }%
    \parindent\z@ \leftskip#3\relax \advance\leftskip\@tempdima\relax
    \rightskip\@pnumwidth plus1em \parfillskip-\@pnumwidth
    #5\leavevmode\hskip-\@tempdima{#6}\nobreak
    \leaders\hbox{$\m@th\mkern \@dotsep mu\hbox{.}\mkern \@dotsep mu$}\hfill
    \nobreak
    \hbox to\@pnumwidth{\@tocpagenum{\ifnum#1=1\bfseries\fi#7}}\par
    \nobreak
    \endgroup
  \fi}
\AtBeginDocument{%
\expandafter\renewcommand\csname r@tocindent0\endcsname{0pt}
}
\def\l@subsection{\@tocline{2}{-5pt}{2.5pc}{5pc}{}}
\makeatother

\title[Approximation of convex functions]{$\mathbf{C^2}$-Lusin approximation of strongly convex functions}

\author[Azagra]{Daniel Azagra}
\address{Daniel Azagra, \newline \indent Department of Mathematical Analysis and Applied Mathematics, \newline \indent Universidad Complutense de Madrid,  28040 Madrid, Spain}
\email{azagra@mat.ucm.es}
\thanks{D.A. was supported by grant PID2022-138758NB-I00.}

\author[Drake]{Marjorie Drake}
\address{Marjorie Drake,\newline \indent Department of Mathematics, Massachusetts Institute of Technology, \newline \indent 77 Massachusetts Ave., Cambridge,
MA 02139}
\email{mkdrake@mit.edu}
\thanks{M.D. was supported by NSF Award No. 2103209.}

\author[Haj\l{}asz]{Piotr Haj\l{}asz}
\address{Piotr Haj\l{}asz,\newline \indent Department of Mathematics, University of Pittsburgh, \newline \indent 301 Thackeray Hall, Pittsburgh,
Pennsylvania 15260}
\email{hajlasz@pitt.edu}
\thanks{P.H. was supported by NSF grant  DMS-2055171 and by Simons Foundation grant 917582.}


\begin{document}

\maketitle
\begin{abstract}
We prove that if $u:\R^n\to\R$ is strongly convex, then for every $\eps>0$ there is a strongly convex function $v\in C^2(\R^n)$ such that $|\{u\neq v\}|<\eps$ and $\Vert u-v\Vert_\infty<\eps$. 
\end{abstract}

\section{Introduction}
\label{intro}

It has been known for at least thirty years, that convex functions have the {\em $C^2$-Lusin property}, meaning that if $u:\R^n\to\R$ is convex, then for every $\eps>0$, there is a function $v\in C^2(\R^n)$ such that $|\{u\neq v\}|<\eps$, see \cite{Alberti2,EvansGangbo,Imomkulov}. Here $|A|$ denotes the Lebesgue measure of $A$. We also say that $v$ is a {\em Lusin approximation} of $u$.

In fact, this result is not particularly difficult. Using Alexandrov's theorem about the second order diffierentiability of $u$, one can show that there is a closed set $E$ contained in the set where $u$ is twice differentiable, with the complement of small measure, $|\R^n\setminus E|<\eps$, and such that the functions $(u|_E,Du|_E,D^2u|_E)$ satisfy the assumptions of the $C^2$-Whitney extension theorem. Since however, the Whitney extension theorem does not preserve convexity, the resulting function $v\in C^2(\R^n)$ that approximates $u$ in the Lusin sense is not convex. 
In fact, in general, one cannot find such convex $v\in C^2(\R^n)$. Here is an easy counterexample:
{\em If $u:\R^2\to\R$, $u(x,y)=|x|$, then the only convex function $v:\R^2\to\R$ such that $|\{u\neq v\}|<\infty$, is $v=u$.}
For this and other examples, see \cite{AH}. In all the examples the problem is caused by the fact that the graph of $u$ contains a line, see \cite[Theorem~1.8]{ACH}, \cite[Proposition~1.10]{AH}.

The problem of approximating a convex function by $C^2$-convex functions in the Lusin sense has remained unresolved since the nineties. It was natural to expect a positive answer in the case of {\em strongly convex functions} or, more generally, {\em locally strongly convex functions}. The motivation for considering this class of functions partly stems from the celebrated work of Greene and Wu \cite{GW}, who proved that locally strongly convex functions on Riemannian manifolds (referred to as {\em strictly convex} by Greene and Wu) can be uniformly approximated by $C^\infty$ locally strongly convex functions. However, it is essential to emphasize that the methods employed by Greene and Wu offer no insight into the problem of Lusin approximation and, in fact, play no role in our paper.

Recall that a function $u: U \to \R$ is {\em strongly convex}, where $U \subseteq \R^n$ is open and convex, if there exists $\eta > 0$ such that $u(x) - \frac{\eta}{2}|x|^2$ is convex. In this case, we say that $u$ is {\em $\eta$-strongly convex}. Moreover, $u$ is {\em locally strongly convex} whenever, for every $x \in U$, there exists $r_x > 0$ such that the restriction of $u$ to $B(x, r_x) \subset U$ is strongly convex. Note that strongly convex functions cannot contain lines on their graphs because they are {\em coercive}, i.e., $u(x) \to \infty$ as $|x| \to \infty$. This excludes the counterexamples discussed above.

The main result of the paper is as follows.
\begin{theorem}
\label{T22}
Let $U\subseteq\R^n$ be open and convex, and $u:U\to\bbbr$ be locally strongly convex. Then for every $\eps_o>0$ and for every continuous function $\varepsilon:U\to (0, 1]$ there is a locally strongly convex function 
$v\in C^2(U)$, such that
\begin{enumerate}
\item[(a)] $|\{x\in U:\, u(x)\neq v(x)\}|<\eps_o$;
\item[(b)] $|u(x)-v(x)|<\varepsilon(x)$ for all $x\in U$;
\item[(c)] $\mathcal{H}^{n}\left( \mathcal{G}_{u} \triangle \mathcal{G}_{v}\right)<\eps_o$.
\end{enumerate}
Also, if $u$ is $\eta$-strongly convex on $U$, then for every $\widetilde{\eta}\in (0, \eta)$ there exists such a function $v$ which is $\widetilde{\eta}$-strongly convex on $U$.
\end{theorem}
Here, $\mathcal{H}^{n}$ denotes Hausdorff $n$-dimensional measure, $\mathcal{G}_{u}$ stands for the graph of a function $u$, and $A\triangle B$ is the symmetric difference of the sets $A$ and $B$, that is, 
$A\triangle B :=(A\setminus B)\cup (B\setminus A)$.

The next result is an immediate consequence of Theorem~\ref{T22}.
\begin{corollary}
\label{main}
Let $u:\bbbr^n\to\bbbr$ be $\eta$-strongly convex. Then for every $\eps>0$ and every $0<\widetilde{\eta}<\eta$, there is an $\widetilde{\eta}$-strongly convex function 
$v\in C^2(\bbbr^n)$, such that
$$
|\{x\in\bbbr^n:\, u(x)\neq v(x)\}|<\eps 
\quad
\text{and}
\quad
\Vert u-v\Vert_\infty<\eps.
$$
\end{corollary}

Let us mention a related result that was proved in \cite[Theorem~1.12]{AH}:
\begin{theorem}
\label{T24}
If $u:\R^n\to\R$ is convex and coercive, then for every $\eps>0$, there is a convex function $v\in C^{1,1}_{\rm loc}(\R^n)$ such that
$|\{x\in\R^n:u(x)\neq v(x)\}|<\eps$.
\end{theorem}
$C^{1,1}_{\rm loc}$ is the class of functions with locally Lipschitz continuous gradient.
In fact due to the counterexamples mentioned above, the coercivity assumption cannot really be removed: it was proved in \cite{AH}, that the existence of Lusin approximation by convex $C^{1,1}_{\rm loc}$ functions as in Theorem~\ref{T24} is equivalent to {\em essential coercivity} of $u$. 
This along with other results in \cite{AH} and \cite[Theorem~1.9]{ACH} provide a complete answer to the problem of Lusin approximation by convex $C^{1,1}$ and $C^{1,1}_{\rm loc}$ functions.

To the best of our knowledge Theorem~\ref{T24} was the first result in the literature towards a positive answer to the question about Lusin approximation by convex functions.

The original proof of Theorem~\ref{T24} was rather difficult and it was based on a characterization of $C_{\rm loc}^{1,1}$-convex Whitney jets \cite[Theorem~1.3]{Azagra2}. However, no results are known for the $C^2$-convex Whitney jets.

A much simpler and more geometric proof of Theorem~\ref{T24} was obtained in \cite{ACH}. There are many geometric conditions that imply $C^{1,1}$ regularity of a convex function. For example it is known that {\em a compact convex body $W$ (i.e. a compact convex set with non-empty interior) has boundary of class $C^{1,1}$ if and only if there is $r>0$ such that $W$ is the union of closed balls of radius $r$}. This result was used in a simple proof of Theorem~\ref{T24} presented in \cite{ACH}. However, we are not aware of any simple geometric conditions that would imply that a convex function is of class $C^2$. 

No methods like those available for $C^{1,1}$-convex functions that were used in \cite{AH,ACH} are available for $C^2$-convex functions. That makes Theorem~\ref{T22} much more difficult. 

In fact our proof of Theorem~\ref{T22} relies on Theorem~\ref{T24} in an essential way. It allows us to assume early in the proof that $u\in C^{1,1}_{\rm loc}$. Then the main work is focused on showing how to upgrade from $C^{1,1}_{\rm loc}$ to $C^2$.
This argument is new and unrelated to methods used in the proof of Theorem~\ref{T24}.

We say that a subset $W$ of $\R^n$ is a {\em locally strongly convex body} if $W$ is closed and convex, with nonempty interior, and $\partial W$ can be locally represented (up to a suitable rotation) as the graph of a strongly convex function. See \cite{Vial} for other equivalent characterizations of locally strongly convex bodies.
\begin{theorem}
\label{geometric corollary}
Let $W \subset \R^n$ be a locally strongly convex body (not necessarily compact), $\varepsilon>0$,  and the set $V \supset \partial W$ be open. There exists a $C^2$ locally strongly convex body $W_{\eps, V}$ such that $\mathcal{H}^{n-1}(\partial W_{\eps, V}\triangle\,\partial W)<\varepsilon$ and $\partial W_{\eps, V}\subset V$. Moreover, if $W$ is compact, then $W_{\eps, V}$ is compact as well.
\end{theorem}
If $W$ is unbounded, then after a suitable rotation, $\partial W$ can be represented as a graph of a locally strongly convex function \cite{AS} and the result follows from Theorem~\ref{T22}. However, if $W$ is compact, the result is more difficult.
It was shown in \cite{ACH,AH} how to prove a related $C^{1,1}$ approximation of convex bodies from Theorem~\ref{T24}. The main idea was to apply Theorem~\ref{T24} to the Minkowski functional of $W$. One can prove that if $W$ is a compact locally strongly convex body, and $\mu_W$ is the Minkowski functional of $W$, then $\mu_W^2$ is strongly convex. Then, a suitable adaptation of methods in \cite{ACH,AH} along with Theorem~\ref{T22} yield Theorem~\ref{geometric corollary}. Details of the proof are, however, quite technical and they will be published elsewhere.

Except for Theorem~\ref{T24} recalled in the Introduction, we use a convention that the names “Theorem” and “Proposition” are reserved for new
results, while well-known results and results of technical character are called “Lemma”. 
“Corollary” will be used for both new and known results, but it will be clear from the context which corollaries are new results.

The paper is structured as follows. In Section~\ref{outline} we outline the main idea of the proof of Theorem~\ref{T22}. In Section~\ref{S2} we collect technical tools that are needed in the proof of Theorem~\ref{T22}.  The results are carefully stated and provide ready to use tools. The reader might browse  though this section quickly and focus on reading 
Section~\ref{S3} that is devoted to 
the proof of Theorem~\ref{T22}. This should give them a motivation for detailed reading of Section~\ref{S2} that otherwise, is a collection of  unrelated results.

\subsection{Notation}
The $L^p$ norm will be denoted by $\Vert u\Vert_{p,E}=(\int_E |u|^p\, dx)^{1/p}$. We will also write $\Vert\cdot\Vert_p$ if the choice of $E$ is obvious. The same convention applies to $p=\infty$ and Sobolev norms discussed below.
$L^p_{\rm loc}(\Omega)$, where $\Omega\subset\bbbr^n$ is open, denotes the space of functions that belong to $L^p(B)$ on every ball $B$ whose closure is contained in $\Omega$, $B\Subset\Omega$. Similar notation applies to Sobolev spaces discussed below.
Lebesgue measure of a set $E$ will be denoted by $|E|$.
The integral average over a ball will be denoted by the barred integral $\mvint_B u\, dx:=|B|^{-1}\int_B u\, dx$. The Lipschitz constant of a function $u$ will be denoted by $\Lip(u)$. The scalar product on $\R^n$ will be denoted by $\langle u,v\rangle$.
By writing $u_n\rightrightarrows u$ we will mean that $u_n$ converges uniformly to $u$.
$E^c:=\bbbr^n\setminus E$ will denote the complement of a set $E$.
By $\nabla^mu$ we will denote the vector whose components are all partial derivatives $D^\alpha u$, $|\alpha|=m$.
$|A|:=\sqrt{\operatorname{tr} (A^TA)}$ will stand for the Hilbert-Schmidt norm of a matrix $A$.
$C$ will denote a generic constant and its value may change from line to line. By writing $C=C(n,m)$ we indicate that the constant depends on parameters $n$ and $m$ only.
$\bbbn$ will denote the set of positive integers.

\section{Outline of the proof of Theorem~\ref{T22}}
\label{outline}

In this section we describe the main ideas in the proof of Theorem~\ref{T22}. 
Throughout the following discussion, we use $\eps$ to represent a small positive constant; its value may vary from one line to another. Since we focus on clarity in this overview, some of the details may be oversimplified.

We begin with the special case that is stated in Lemma~\ref{T25}: {\em If $u:\R^n\to\R$ is strongly convex, then there is a strongly convex $v\in C^2(\R^n)$ such that $|\{ u\neq v\}|<\eps$.}

As a consequence of Theorem~\ref{T24}, we may assume that $u\in C^{1,1}_{\rm loc}(\R^n)$, and hence, $u$ belongs to any Sobolev spaces $W^{2,p}_{\rm loc}(\R^n)$ for $1\leq p\leq\infty$. In Proposition~\ref{T6}, we decompose $\R^n$ into annuli where we approximate $u$ by its convolution; we patch these approximations together with a smooth partition of unity to produce a strongly convex function $\widetilde{u}\in C^\infty(\R^n)$ satisfying $u-\widetilde{u}\in W^{2,1}(\R^n)$ and $\Vert u-\widetilde{u}\Vert_{2,1}<\eps$. 

The function $\widetilde{u}$ need not agree with $u$, but we can use the small magnitude of $\Vert u-\widetilde{u}\Vert_{2,1}$ to prove there exists $w\in C^2(\R^n)$ with small seminorm $\Vert\nabla^2 w\Vert_\infty$ satisfying $|\{w\neq u-\widetilde{u}\}|<\eps$. Then we have $v:=w+\widetilde{u}\in C^2$, a strongly convex function by Lemma~\ref{T12}, and $|\{ v\neq u\}|=|\{w\neq u-\widetilde{u}\}|<\eps$. Finding $w$ is tricky:

First, let $\hat{u}:=u-\widetilde{u}\in W^{2,1}(\R^n)$, so $\Vert\hat{u}\Vert_{2,1}<\eps$. Pointwise inequalities from \cite{BH,BHS} (see Corollary~\ref{T2}) imply that
$$
\frac{|\hat{u}(y)-\hat{u}(x)-D\hat{u}(x)(y-x)|}{|y-x|^2}+\frac{|D\hat{u}(y)-D\hat{u}(x)|}{|y-x|}\leq C\big(\MM|\nabla^2 \hat{u}|(x)+\MM|\nabla^2\hat{u}|(y)\big),
$$
where $\MM g$ is the Hardy-Littlewood maximal function. 
We choose $\lambda>0$ large enough that 
$$
E_\lambda:=\{\MM|\nabla^2\hat{u}|\leq\lambda\}
\quad
\text{satisfies}
\quad
|\R^n\setminus E_\lambda|<\eps.
$$
Then $\hat{u}$ restricted to $E_\lambda$ satisfies the assumptions of the $W^{2,\infty}$-Whitney extension theorem (Lemma~\ref{T3}); thus, there exists $\hat{w}\in W^{2,\infty}(\R^n)\subset C^{1,1}(\R^n)$ satisfying $\hat{w}|_{E_\lambda}=\hat{u}|_{E_\lambda}$, and 
$$
|\{\hat{w}\neq \hat{u}\}|\leq |\R^n\setminus E_\lambda|<\eps.
$$ 
We let $\hat{v} :=\hat{w}+\tilde{u}\in C^{1,1}_{\rm loc}$. Even though $\hat{v}$ has the same regularity as $u$, we have $|\{\hat{v} \neq u\}|=|\{\hat{w}\neq u-\widetilde{u}\}|<\eps$. Further, we are able to use our control of $\|\hat{u} \|_{2,1}$, to control $\|\nabla^2\hat{w}\|_{\infty}$, and in a delicate iterative procedure (carefully explained in the proof of Proposition~\ref{T4}),  
we produce a series of $C^\infty$ functions that converge in the $W^{2,\infty}$ norm to a function $w$ that has small  seminorm $\Vert\nabla^2 w\Vert_\infty$ and satisfies that $|\{w\neq u-\widetilde{u}\}|<\eps$. In particular, the function $w\in C^2$, because functions in the closure of $C^\infty$ in the $W^{2,\infty}$ norm are of class $C^2$. By letting $v:= w+\widetilde{u}$ be our approximation for $u$, we complete the proof of Lemma~\ref{T25}. 

In Theorem~\ref{T22}, we extend this result to approximation of {\em locally} strongly convex functions on a convex domain $U$, where in addition to the Lusin property (a) we have the uniform approximation (b) ((c) follows easily from (a) and the local Lipschitz continuity of convex functions). This adds another layer of difficulty. We represent $U$ as the union of an increasing sequence of convex bodies
$$
B_k\subset\inte(C_k)\subset C_k\subset\inte(B_{k+1}),
\qquad
\bigcup_{k=1}^\infty B_k=U.
$$
Using Corollary~\ref{T18}, we extend $u$ from $B_k\setminus C_{k-1}$ to a strongly convex function $u_k:\R^n\to\R$. Then, we approximate each $u_k$ by a strongly convex function $h_k\in C^2(\R^n)$ using Lemma~\ref{T25}. We glue the functions $h_k$ using the smooth maxima described in Lemma~\ref{T17}, so that the resulting function $v\in C^2(\Omega)$ equals $h_k$ on $B_k\setminus C_{k-1}$. This function $v$ satisfies (a)-(c) of Theorem~\ref{T22}.

\section{Preliminaries for the proof of Theorem~\ref{T22}}
\label{S2}

\subsection{Smooth and Sobolev functions}
Let $\Omega\subset\bbbr^n$ be open, $m\in\bbbn$, and $1\leq p\leq\infty$.

$C^m(\Omega)$ is the space of functions on $\Omega$ that are continuously $m$-times differentiable and $C_0^\infty(\Omega)$ is the space of compactly supported smooth functions.
$C^{m,1}(\Omega)$ is the subspace of $C^m(\Omega)$ consisting of functions $u$ such that the derivatives $D^\alpha u$, $|\alpha|=m$, are Lipschitz continuous on $\Omega$, and $C^{m,1}_{\rm loc}(\Omega)$ is the space of functions whose derivatives of order $m$ are locally Lipschitz continuous in $\Omega$.

We will use only standard facts from the theory of Sobolev spaces. For more details and proofs of results listed here, see e.g. \cite{EG}.

$W^{m,p}(\Omega)$
is the Sobolev spaces of functions whose weak derivatives satisfy
$$
\Vert u\Vert_{m,p,\Omega}:=\sum_{|\alpha|\leq m}\Vert D^\alpha u\Vert_{p,\Omega}<\infty.
$$
As explained above, for a measurable set $E\subset\Omega$ we will write
$$
\Vert u\Vert_{m,p,E}:=\sum_{|\alpha|\leq m}\Vert D^\alpha u\Vert_{p,E}.
$$
However, $\Vert D^2 u\Vert_p$ will be the $L^p$ norm of $|D^2u|$ i.e., 
$
\Vert D^2 u\Vert_p = \Vert(\sum_{|\alpha|=2}|D^\alpha u|^2)^{1/2}\Vert_p.
$
Clearly, $\Vert D^2 u\Vert_p\leq \Vert u\Vert_{2,p}$. 

If $\Omega\subset\bbbr^n$ is convex, then $W^{m,\infty}(\Omega)\subset C^{m-1,1}(\Omega)$. Indeed, convexity of $\Omega$ implies that $W^{1,\infty}(\Omega)$ coincides with the space of bounded Lipschitz continuous functions on $\Omega$ and hence if $u\in W^{m,\infty}(\Omega)$, and $|\alpha|=m-1$, then $D^\alpha u\in W^{1,\infty}(\Omega)$ is Lipschitz continuous.

Note that $W^{m,\infty}(\Omega)\cap C^m(\Omega)$ is a closed subspace of $W^{m,\infty}(\Omega)$. Indeed, if $\{ u_k\}_k\subset W^{m,\infty}\cap C^m(\Omega)$ is a Cauchy sequence with respect to the norm $\Vert\cdot\Vert_{m,\infty}$, all derivatives $D^\alpha u_k$, $|\alpha|\leq m$ converge uniformly on $\Omega$ and hence the limiting function belongs to $C^m(\Omega)$.

We also define $W^m(\Omega):=W^{m,1}(\Omega)\cap W^{m,\infty}(\Omega)$. Thus $u\in W^m(\Omega)$ if $D^\alpha u\in L^p(\Omega)$ for all $1\leq p\leq\infty$ and all $|\alpha|\leq m$. Indeed, this is a direct consequence of an elementary fact that $L^1\cap L^\infty\subset L^p$ for all $1\leq p\leq\infty$. 

\subsection{Mollifiers}
\label{S1}
By a mollifier we will mean a compactly supported smooth function $\vi\in C_0^\infty(B(0,1))$ such that $\vi\geq 0$ and $\int_{\bbbr^n}\vi\, dx=1$. Then for $\eps>0$ we set $\vi_\eps(x)=\eps^{-n}\vi(x/\eps)$. The following facts are well known and easy to prove.

If $u\in W^{m,p}(\bbbr^n)$, then $u_\eps:=\vi_\eps*u\in C^\infty(\bbbr^n)$ satisfy
\begin{equation}
\label{eq14}
\Vert u-u_\eps\Vert_{m,p}\stackrel{\eps\to 0^+}{\longrightarrow} 0
\quad
\text{if $1\leq p<\infty$}
\end{equation}
and
\begin{equation}
\label{eq15}
\Vert u_\eps\Vert_{m,p}\leq \Vert u\Vert_{m,p}
\quad
\text{if $1\leq p\leq \infty$.}
\end{equation}
It is also easy to show that $\Vert D^2 u_\eps\Vert_\infty\leq \Vert D^2 u\Vert_\infty$.
\subsection{Convex functions}
For a function $u:U\to\R$ defined on an open and convex set $U\subseteq\R^n$ and $x\in U$, we define the {\em subdifferential} $\partial u(x)$ as the set of all $v\in\R^n$ such that $u(y)\geq u(x)+\langle v,y-x\rangle$ for all $y\in U$. 
If $A\subset U$, then we write $\partial u(A)=\bigcup_{x\in A}\partial u(x)$.

It is well known that $u$ is convex if and only if $\partial u(x)\neq\varnothing$ for all $x\in U$. 

The next result is an easy exercise.
\begin{lemma}
\label{T23}
Let $u:U\to\R$ be a convex function defined on an open and convex set $U\subseteq\R^n$. If $G\Subset U$ is open, then $\sup_{\xi\in\partial u(G)}|\xi|\leq\Lip(u|_G)$.
\end{lemma}
We will also need the following results.
\begin{lemma}
\label{T14}
Let $U\subseteq\R^n$ be open and convex, and $u_i:U\to\R$ be a sequence of convex functions. Assume that the sequence $(u_i)$ converges pointwise on a dense subset of $U$. Then there is a convex function $u:U\to\R$ such that $u_i\rightrightarrows u$ converges uniformly on every compact subset of $U$.
\end{lemma}
For a proof, see \cite[Theorem 10.8]{Rockafellar}.
\begin{lemma}
\label{T13}
Let $U\subseteq \R^n$ be open and convex, $u:U\to\R^n$ be a convex function, and $A, B \subset U$ be compact convex bodies such that $A\subset\inte(B)$. Then, 
\begin{equation}
\label{eq56}
\Lip(u_{|A})\leq \frac{M_B(u) -m_B(u)}{\dist(\partial A,\partial B)},
\end{equation} 
where $M_B(u):=\max_{x\in B}u(x)$, and $m_B(u):=\min_{x\in B}u(x)$.
\end{lemma}
\begin{proof}
Let $x,y\in A$, $x\neq y$. Without loss of generality, we may assume that $u(y)\geq u(x)$. Then, there exists $z\in\partial B$ such that $y$ belongs to the segment with the endpoints $x$ and $z$, and we have
$$
\frac{u(y) - u(x)}{|y-x|} \leq \frac{u(z)-u(x)}{|z-x|} \leq \frac{M_B(u) -m_B(u)}{\dist(\partial A,\partial B)},
$$
implying \eqref{eq56}.
\end{proof}

\subsection {Strongly convex functions}
Let $U\subseteq\R^n$ be open and convex. 
We say that a function $u:U\to\bbbr$ is {\em $\eta$-strongly convex}, where $\eta>0$ is a constant, if the function $u(x)-\frac{\eta}{2}|x|^2$ is convex. A function is {\em strongly convex} if it is $\eta$-strongly convex for some $\eta>0$.
\begin{lemma}
\label{T16}
Let $u:U\to\R$ be a convex function defined on an open and convex set $U\subseteq\R^n$. Then $u$ is $\eta$-strongly convex if and only if 
\begin{equation}
\label{eq39}
u(y)\geq u(x)+\langle\xi,y-x\rangle +\frac{\eta}{2}|y-x|^2
\end{equation}
for all $x,y\in U$ and $\xi\in \partial u(x)$.
\end{lemma}
\begin{proof}
Let $v(x)=u(x)-\frac{\eta}{2}|x|^2$ and note that \eqref{eq39} is equivalent to 
\begin{equation}
\label{eq40}
v(y)\geq v(x)+\langle \xi-\eta x,y-x\rangle.
\end{equation}
To prove the implication $\Leftarrow$, let $u:U\to\R$ be convex and satisfy \eqref{eq39} or equivalently \eqref{eq40}.
This implies that $\xi-\eta x\in \partial v(x)$ so $\partial v(x)\neq\varnothing$ for all $x\in U$, and hence $v$ is convex i.e., $u$ is $\eta$-strongly convex.
It remains to prove the implication $\Rightarrow$. Let $u:U\to\R$ be $\eta$-strongly convex, let $x,y\in U$ and let $\xi\in \partial u(x)$. We need to prove \eqref{eq39} or equivalently, \eqref{eq40}. Let $z\in U$. Subtracting the identity
$$
\frac{\eta}{2}|z|^2=\frac{\eta}{2}|x|^2+\langle\eta x,z-x\rangle +\frac{\eta}{2}|z-x|^2
$$
from the inequality
$u(z)\geq u(x)+\langle\xi,z-x\rangle$, yields
\begin{equation}
\label{eq41}
v(z)\geq v(x)+\langle\xi-\eta x,z-x\rangle -\frac{\eta}{2}|z-x|^2.
\end{equation}
The function
$$
h(z):=v(z)-\langle\xi-\eta x,z-x\rangle
$$
is convex and \eqref{eq41} is equivalent to
\begin{equation}
\label{eq42}
h(z)\geq h(x)-\frac{\eta}{2}|z-x|^2.
\end{equation}
For $t\in (0,1)$, let $z=x+t(y-x)$. Then \eqref{eq42} and convexity of $h$ yield
$$
th(y)+(1-t)h(x)\geq h(x+t(y-x))\geq h(x)-\frac{\eta}{2}|t(y-x)|^2
$$
i.e., $h(x)-h(y)\leq t\frac{\eta}{2}|y-x|^2$. Letting $t\to 0^+$ gives $h(y)\geq h(x)$, which is precisely inequality \eqref{eq40}.
\end{proof}
\begin{remark}\label{remark on characterization of strong convexity}
The above proof of $(\Leftarrow)$ also shows that if \eqref{eq39} holds for all $x, y\in U$ and {\em some } $\xi\in \partial u(x)$, then $u$ is $\eta$-strongly convex, and therefore, by the proof of $(\Rightarrow)$, \eqref{eq39} is also true for {\em all } $\xi\in\partial u(x)$.
\end{remark}

As a corollary we obtain a result about the minimal $\eta$-strongly convex extension from a compact set.
\begin{corollary}
\label{T18}
Let $u:U\to\R$ be an $\eta$-strongly convex function defined on an open and convex set $U\subseteq R^n$. If $K\subset U$ is compact, then the function $\widetilde{u}:\R^n\to\R$ defined by
\begin{equation}
\label{eq43}
\widetilde{u}(x):=\sup_{z\in K,\xi_z\in \partial u(z)} 
\Big\{ u(z)+\langle\xi_z,x-z\rangle+\frac{\eta}{2}|x-z|^2\Big\}
\end{equation}
is $\eta$-strongly convex, and such that
\begin{enumerate}
\item[(a)] $\widetilde{u}(x)=u(x)$ for all $x\in K$,
\item[(b)] $\widetilde{u}(x)\leq u(x)$ for all $x\in U$.
\end{enumerate}
\end{corollary}
\begin{proof}
It follows from Lemma~\ref{T23} and Lemma~\ref{T13} that $\partial u(K)$
is bounded. This implies that $\widetilde{u}(x)$ is finite for every $x\in\R^n$. The function $\widetilde{u}(x)-\frac{\eta}{2}|x|^2$ is convex as the supremum of a collection of affine functions so $\widetilde{u}$ is $\eta$-strongly convex. Finally, Lemma~\ref{T16} implies that $u(x)\geq\widetilde{u}(x)$ for all $x\in U$ which is (b), and taking $z=x\in K$ and any $\xi_z\in\partial u(z)$ in \eqref{eq43} yields $\widetilde{u}(x)\geq u(x)$ for $x\in K$ so $\widetilde{u}(x)=u(x)$ for all $x\in K$ which is (a). 
\end{proof}

\subsection{Smooth strongly convex functions}
Note that if $u\in C^2(U)$, then $u$ is $\eta$-strongly convex if and only if
$$
D^2u(x)\geq\eta
\qquad
\text{for all $x\in\bbbr^n$,}
$$
meaning that
\begin{equation}
\label{eq32}
\xi^TD^2u(x)\xi :=
\sum_{i,j=1}^n\frac{\partial^2 u}{\partial x_i\partial x_j}(x)\xi_i\xi_j\geq\eta|\xi|^2
\quad
\text{for all $x,\xi\in\bbbr^n$.}
\end{equation}
\begin{lemma}
\label{T5}
Let $u:\bbbr^n\to\bbbr$ be convex and $u_\eps:=u*\vi_\eps$, where $\vi$ is a mollifier and $\eps>0$. Then
\begin{enumerate}
\item[(a)] $u_\eps$ is convex,
\item[(b)] If $u$ is $\eta$-strongly convex, then $u_\eps$ is $\eta$-strongly convex too.
\end{enumerate}
\end{lemma}
\begin{proof}
(a) Multiplying the inequality 
$u((tx+(1-t)y)-z)\leq tu(x-z)+(1-t)u(y-z)$,
$t\in [0,1]$,
by $\vi_\eps(z)$ and integrating over $\bbbr^n$ yields
$u_\eps(tx+(1-t)y)\leq tu_\eps(x)+(1-t)u_\eps(y)$.

\noindent (b)
Suppose first that $u\in C^\infty(\bbbr^n)$ is $\eta$-strongly convex. Then \eqref{eq32} is satisfied and hence
$$
\xi^TD^2u_\eps \xi=(\xi^TD^2u\, \xi)*\vi_\eps\geq \eta|\xi|^2.
$$
This proves $\eta$-strong convexity of $u_\eps$ in the smooth case.

In particular, if $g(x)=\frac{\eta}{2}|x|^2$, then
$$
\xi^TD^2g_\eps \xi\geq \eta|\xi|^2.
$$
Assume now that $u$ is any $\eta$-strongly convex function. Then $(u-g)_\eps=u_\eps-g_\eps$ is convex by (a) and hence 
$$
\xi^TD^2(u_\eps-g_\eps)\xi\geq 0,
\qquad
\xi^TD^2u_\eps \xi\geq \xi^TD^2g_\eps \xi\geq\eta|\xi|^2.
$$
This proves $\eta$-strong convexity of $u_\eps$ in the general case.
\end{proof}
\begin{lemma}
\label{T12}
If $u:\R^n\to\R$ is $\eta$-strongly convex, $0<\widetilde{\eta}<\eta$, and $f\in C^2(\bbbr^n)$ satisfies
$$
n\max_{|\alpha|=2}\Vert D^\alpha f\Vert_\infty\leq\eta-\widetilde{\eta},
$$
then $u+f$ is $\widetilde{\eta}$-strongly convex.
\end{lemma}
\begin{proof}
The following elementary estimate
$$
|v^T D^2f\, v|\leq (\eta-\widetilde{\eta})|v|^2,
\quad
v\in\bbbr^n,
$$
implies that the function $f(x)+\frac{\eta-\widetilde{\eta}}{2}|x|^2$ is convex. Therefore,
$$
(u(x)+f(x))-\frac{\widetilde{\eta}}{2}|x|^2=\Big(u(x)-\frac{\eta}{2}|x|^2\Big)+\Big(f(x)+\frac{\eta-\widetilde{\eta}}{2}|x|^2\Big)
$$
is convex.
\end{proof}

\subsection{$C^{1,1}$ approximation of strongly convex functions}
We say that a function $u:\bbbr^n\to\bbbr$ is
{\em coercive} if 
$$
\lim_{|x|\to\infty} u(x)=\infty. 
$$
Note that strongly convex functions are coercive. Indeed, convex functions are bounded from below by affine functions (supporting hyperplanes) and hence strongly convex functions are bounded from below by quadratic functions.

The following result is a straightforward consequence of \cite[Theorem~1.12]{AH}. For a geometric and a much simpler proof, see \cite{ACH}.
\begin{lemma}
\label{T10}
If $u:\bbbr^n\to\bbbr$ is convex and coercive, then for every $\eps>0$, there is a convex function $v\in C^{1,1}_{\rm loc}(\bbbr^n)$ such that\footnote{Although the proof of \cite{AH} implicitly yields $v\geq u$, this useful additional property is explicitly stated and proved in \cite{ACH}.} $|\{x\in\bbbr^n:\, u(x)\neq v(x)\}|<\eps$.
\end{lemma}
\begin{corollary}
\label{T11}
If $u:\bbbr^n\to\bbbr$ is $\eta$-strongly convex, then for every $0<\widetilde{\eta}<\eta$ and every $\eps>0$, there is a 
$\widetilde{\eta}$-strongly convex function $v\in C^{1,1}_{\rm loc}(\bbbr^n)$, such that $v\geq u$ and $|\{x\in\bbbr^n:\, u(x)\neq v(x)\}|<\eps$.
\end{corollary}
\begin{proof}
The function $\widetilde{u}(x)=u(x)-(\widetilde{\eta}/2)|x|^2$ is strongly convex and hence coercive. Applying
Lemma~\ref{T10} to that function yields a convex function $\widetilde{v}\in C^{1,1}_{\rm loc}$ satisfying $\widetilde{v}\geq\widetilde{u}$ and $|\{\widetilde{u}\neq\widetilde{v}\}|<\eps$. Therefore, $v=\widetilde{v}+(\widetilde{\eta}/2)|x|^2$ satisfies the claim of the corollary.
\end{proof}
\subsection{Smooth maximum}
The material of this subsection is based on \cite{Azagra}. The maximum of two convex functions $\max\{u,v\}=(u+v+|u-v|)/2$ is convex, but not necessarily differentiable, even if $u,v\in C^\infty$. The next construction provides a modified notion of the maximum that preserves smoothness and convexity.

For any constant $\delta>0$, we define
$$
M_\delta(x,y):=\frac{x+y+\theta(x-y)}{2},
\quad(x,y)\in \R^2,
$$
where $\theta=\theta_{\delta}\colon\R\to (0,
\infty)$ is a $C^\infty$ function such that $\theta(t)=|t|$ if and only if $|t|\geq\delta$,
and $\theta$ is convex and symmetric. It follows that $M_\delta:\R^2\to\R$ is convex.

For functions $u, v\colon U\subseteq\R^n\to\R$, we define the {\em smooth maximum} of $u$ and $v$  by
\begin{equation}
\label{eq57}
M_\delta(u,v):U\to\R,
\qquad
M_{\delta}(u,v)(x)=M_{\delta}(u(x), v(x)).
\end{equation}
We will need the following properties of the smooth maximum. For this and other properties, see \cite[Section~2]{Azagra}.

\begin{lemma}
\label{T17}
Let $U\subseteq\R^n$ be open and convex and let $u, v\colon U\to\R$
be convex functions. For every $\delta>0$, the function
$M_{\delta}(u,v)\colon U\to\R$ is convex and satisfies
\begin{enumerate}
\item[(a)] If $u,v\in C^2(U)$, then $M_{\delta}(u,v)\in C^2(U)$.
\item[(b)] $M_{\delta}(u,v)(x)=u(x)$ if $u(x)\geq v(x)+\delta$, and $M_{\delta}(u,v)(x)=v(x)$ if $v(x)\geq u(x)+\delta$.
\item[(c)] $\max\{u,v\}\leq M_{\delta}(u,v)\leq \max\{u,v\} + \delta/2$.
\item[(d)] If $u,v$ are $\eta$-strongly convex for some $\eta>0$, then $M_{\delta}(u,v)$ is $\eta$-strongly convex.
\end{enumerate}
\end{lemma}
\begin{proof}
It is not difficult to show that the partial derivatives of $M_\delta(x,y)$ are non-negative and hence $M_\delta(x,y)$ is non-decreasing in $x$ and in $y$, see \cite[Lemma~2(viii)]{Azagra}. This easily implies convexity of $M_\delta(u,v)$. Properties (a)-(c) follow immediately from the definition of $M_\delta(x,y)$. It remains to prove property (d). If $u,v$ are $\eta$-strongly convex, and $h(x) = (\eta/2) |x|^2$, then $M_\delta(u,v) - h=M_\delta(u-h,v-h)$ is convex, so $M_\delta(u,v)$ is $\eta$-strongly convex. 
\end{proof}

\subsection{Smooth approximation of strongly convex functions}
\begin{proposition} 
\label{T6}
Let $u \in W^{2,p}_{\loc}(\bbbr^n)$, $n<p<\infty$, be $\eta$-strongly convex. Then for every $\eps>0$, $0<\widetilde{\eta}<\eta$, and $1\leq q\leq p$, there exists an $\widetilde{\eta}$-strongly convex function $v \in C^\infty(\bbbr^n)$ such that
$\Vert u - v\Vert_{2,q,\bbbr^n} < \eps$.
\end{proposition}
\begin{proof}
Let $\psi_i\in C_0^\infty\big((i-1,i+1)\big)$, $i=0,1,2,\ldots$ be such that
$\sum_{i=0}^\infty\psi_i(t)=1$ for all $t\in [0,\infty)$. Then the functions $\theta_i\in C_0^\infty(\bbbr^n)$, $\theta_i(x)=\psi_i(|x|)$ form a partition of unity on $\bbbr^n$. The function $\theta_0$ is supported in $B(0,1)$, while the functions $\theta_i$ for $i\geq 1$ are supported in the annuli $\{i-1< |x|<i+1\}$.

Let $A_0=\overbar{B}(0,1)$, and $A_i=\{i-1\leq |x|\leq i+1\}$ for $i\geq 1$, so $\supp\theta_i\subset A_i$.

Fix a mollifier $\vi$ and let $u_\delta:=u*\vi_\delta$.
Note that for every $i$, $\Vert u-u_\delta\Vert_{2,p,A_i}\to 0$ as $\delta\to 0^+$, and hence
\begin{equation}
\label{eq33}
\Vert u-u_\delta\Vert_{2,q,A_i}\to 0
\qquad
\text{as $\delta\to 0$,}
\end{equation}
by H\"older's inequality.
Since $u$ is continuous and $Du\in W^{1,p}_{\rm loc}$, $p>n$, is continuous by the Sobolev embedding theorem, we have that for every $i$,
\begin{equation}
\label{eq35}
\Vert u- u_\delta\Vert_{\infty, A_i}+\Vert D(u-u_\delta)\Vert_{\infty, A_i}\to 0
\qquad
\text{as $\delta\to 0^+$.}
\end{equation}

Fix $\eps>0$, $0<\widetilde{\eta}<\eta$ and $1\leq q\leq p$.

Let $\Delta=\{\delta_i\}_{i=0}^\infty$ be a sequence of positive numbers $\delta_i>0$, and define
$$
u_\Delta:=\sum_{i=0}^\infty\theta_i u_{\delta_i}.
$$
We will prove that with a suitable choice of $\Delta$, the function $v:=u_\Delta$ has the required properties. 

We say that $\Delta\leq\Delta'=\{\delta'_i\}_{i=0}^\infty$, if $\delta_i\leq \delta_i'$ for all $i$. We say that a property holds for all sufficiently small $\Delta$, if there is $\Delta'$ such that the property holds for all $\Delta\leq\Delta'$.

Clearly $u_\Delta\in C^\infty(\bbbr^n)$ for every $\Delta$. Next, we show that there is $\Delta'$ such that
\begin{equation}
\label{eq34}
\Vert u-u_\Delta\Vert_{2,q,\bbbr^n}< \eps
\qquad
\text{for all $\Delta\leq\Delta'$.}
\end{equation}
Using \eqref{eq33} we can conclude that
$$
\Vert u-u_\Delta\Vert_{q,\bbbr^n}\leq \sum_{i=0}^\infty \Vert\theta_i(u-u_{\delta_i})\Vert_{q,\bbbr^n}\leq
\sum_{i=0}^\infty \Vert u-u_{\delta_i}\Vert_{q,A_i}<\eps/3,
$$
provided $\Delta$ is sufficiently small. The estimate for the gradient is similar:
$$
\Vert D(u-u_\Delta)\Vert_{q,\bbbr^n}\leq \sum_{i=0}^\infty \Vert D\theta_i\Vert_\infty \Vert u-u_{\delta_i}\Vert_{q,A_i}+
\sum_{i=0}^\infty\Vert D(u-u_{\delta_i})\Vert_{q,A_i}<\eps/3,
$$
provided $\Delta$ is sufficiently small.

To estimate the second order derivatives we will use the following convenient notation for the product rule. If $a=[a_i]_{i=1}^n$ and $b=[b_i]_{i=1}^n$, then the tensor product of the vectors is the $n\times n$ matrix defined by $a\otimes b:=[a_ib_j]_{i,j=1}^n$. With this notation, the product rule for the second order derivatives on $\bbbr^n$ takes the form of
$$
D^2(fg)=gD^2f+Df\otimes Dg+ Dg\otimes Df+fD^2g.
$$
Thus the Hilbert-Schmidt norm of $D^2(fg)$ satisfies
$$
|D^2(fg)|\leq |g|\, |D^2f|+2|Df|\, |Dg|+|f|\, |D^2g|
$$
and hence,
\begin{equation*}
\begin{split}
&\Vert D^2(u-u_\Delta)\Vert_{q,\bbbr^n} \\
&\leq 
\sum_{i=0}^\infty\Big(\Vert D^2(u-u_{\delta_i})\Vert_{q,A_i}+2\Vert D\theta_i \Vert_\infty\Vert D(u-u_{\delta_i})\Vert_{q,A_i}+
\Vert D^2\theta_i\Vert_\infty \Vert u-u_{\delta_i}\Vert_{q,A_i}\Big)<\frac{\eps}{3},
\end{split}
\end{equation*}
provided $\Delta$ is sufficiently small. This proves \eqref{eq34} and it remains to show that $u_\Delta$ is $\widetilde{\eta}$-strongly convex for all sufficiently small $\Delta$. 

Since $\sum_{i=0}^\infty \theta_i=1$, we have that $\sum_{i=0}^\infty D\theta_i=0$, $\sum_{i=0}^\infty D^2\theta_i=0$ and hence
$$
\sum_{i=0}^\infty \big(D\theta_i\otimes Du + Du \otimes D\theta_i+  uD^2\theta_i\big)=0.
$$
Therefore,
$$
D^2u_\Delta=
\underbrace{\sum_{i=0}^\infty \theta_i D^2u_{\delta_i}}_{A_\Delta}\, +\, 
\underbrace{\sum_{i=0}^\infty D\theta_i\otimes D(u_{\delta_i}-u)+D(u_{\delta_i}-u)\otimes D\theta_i}_{B_\Delta}\, +\, 
\underbrace{\sum_{i=0}^\infty(u_{\delta_i}-u)D^2\theta_i}_{C_\Delta}
$$
According to Lemma~\ref{T5}, the functions $u_{\delta_i}$ are $\eta$-strongly convex and hence 
$$
\xi^T A_\Delta \xi\geq \eta |\xi|^2,
\qquad
\xi\in\bbbr^n.
$$
On the other hand \eqref{eq35} yields that for sufficiently small $\Delta$ we have
$$
|\xi^T B_\Delta\, \xi|\leq 2|\xi|^2\sum_{i=0}^\infty\Vert D\theta_i\Vert_\infty \Vert D(u_{\delta_i}-u)\Vert_{\infty,A_i}\leq
\frac{\eta-\widetilde{\eta}}{2}\, |\xi|^2,
$$
and 
$$
|\xi^T C_\Delta\, \xi|\leq |\xi|^2\sum_{i=0}^\infty \Vert D^2\theta_i\Vert_\infty\Vert u_{\delta_i}-u\Vert_{\infty,A_i}\leq
\frac{\eta-\widetilde{\eta}}{2}\, |\xi|^2.
$$
In concert, the last three estimates yield
$$
\xi^T D^2 u_\Delta\, \xi\geq \eta|\xi|^2-\frac{\eta-\widetilde{\eta}}{2}\, |\xi|^2-\frac{\eta-\widetilde{\eta}}{2}\, |\xi|^2= 
\widetilde{\eta}|\xi|^2,
$$
and that completes the proof of $\widetilde{\eta}$-strong convexity of $u_\Delta$.
\end{proof}

\subsection{Pointwise inequality}
Given $g\in L^1_{\rm loc}(\bbbr^n)$, the Hardy-Littlewood maximal function is defined on $\bbbr^n$ by
$$
\MM g(x):= \sup_{r>0} \mvint_{B(x,r)} |g(y)|\, dy.
$$
Recall that the maximal function satisfies the weak type estimate \cite[p. 5]{stein}:
\begin{equation}
\label{eq11}
|\{x\in\bbbr^n:\, \MM g(x)>t\}|\leq\frac{5^n}{t}\Vert g\Vert_1
\quad
\text{for all $t>0$.}
\end{equation}

Sobolev functions $u\in W^{m,1}_{\rm loc}(\Omega)$ and their derivatives are defined almost everywhere. In what follows we will choose representatives defined {\em everywhere} in $\Omega$ by the following formula:
\begin{equation}
\label{eq1}
D^\alpha u(x):=\limsup_{r\to 0}\mvint_{B(x,r)} D^\alpha u(y)\, dy,
\quad
\text{ for all $|\alpha|\leq m$ and every $x\in\Omega$.}
\end{equation}
The fact that it is indeed a representative of $D^\alpha u$ in the class of functions equal a.e. follows immediately from the Lebesgue differentiation theorem.

For $k\leq m$, we will denote the Taylor polynomial of $u\in W^{m,1}_{\rm loc}(\Omega)$
at $x\in\Omega$ by
$$
T_x^ku(y):=\sum_{|\alpha|\leq k} D^\alpha u(x)\frac{(y-x)^\alpha}{\alpha!}\, .
$$
More generally, if $\{ u^\alpha\}_{|\alpha|\leq m}$ is a family of continuous functions on $F\subset\bbbr^n$, then for $|\alpha|\leq m$ we define
$$
T_x^ku^\alpha(y):= \sum_{|\beta|\leq k} u^{\alpha+\beta}(x)\frac{(y-x)^\beta}{\beta !}
\quad
\text{for $x\in F$, $y\in\bbbr^n$ and $k\leq m-|\alpha|$.}
$$

For a proof of the following pointwise inequality, see \cite[Theorem~3]{BH} and \cite[Theorem~3.6]{BHS}.
In the case of smooth functions this result has already been proved in \cite[Corollary~5.8]{BS}.
\begin{lemma}
\label{T1}
Let $u\in W^{m,1}_{\rm loc}(\bbbr^n)$ and its all derivatives  be defined at every point of $\bbbr^n$ by \eqref{eq1}. Then
\begin{equation}
\label{eq2}
\frac{|u(y)-T_x^{m-1}u(y)|}{|x-y|^{m}}\leq C(\MM |\nabla^mu|(x)+\MM|\nabla^mu|(y))
\end{equation}
for all $x,y\in\bbbr^n$, $x\neq y$, where the constant $C$ depends on $n$ and $m$ only.
\end{lemma}
\begin{remark}
If $u(y)=\pm\infty$ or $D^\alpha u(x)=\pm\infty$ for some $|\alpha|\leq m-1$, we assume that the left hand side of \eqref{eq2} equals $+\infty$ so the inequality means that in that case the right hand side equals $+\infty$. 
\end{remark}
Applying Lemma~\ref{T1} to $D^\alpha u\in W^{m-|\alpha|,1}_{\rm loc}$, $|\alpha|\leq m-1$, we have
\begin{corollary}
\label{T2}
Let $u\in W^{m,1}_{\rm loc}(\bbbr^n)$ and its all derivatives  be defined at every point of $\bbbr^n$ by \eqref{eq1}. Then we have
\begin{equation}
\label{eq3}
\max_{|\alpha|\leq m-1}\frac{|D^\alpha u(y)-T_x^{m-1-|\alpha|}D^\alpha u(y)|}{|x-y|^{m-|\alpha|}}\leq C_p(\MM |\nabla^mu|(x)+\MM|\nabla^mu|(y))
\end{equation}
for all $x,y\in\bbbr^n$, $x\neq y$, where the constant $C_p$ depends on $n$ and $m$ only.
\end{corollary}

Lemma~\ref{T1} and Corollary~\ref{T2} provide estimates for the Taylor reminders of $u$ and its derivatives in terms of the maximal function. These estimates allow us to apply the Whitney extension theorem which we formulate next. A proof of this particular statement can be found in \cite[p. 177--180]{stein}.
\begin{lemma}
\label{T3}
Let $m\in\bbbn$ and $M>0$. Assume that $F\subset\bbbr^n$ is a closed set and we have a collection of continuous functions $\{u^\alpha\}_{|\alpha|\leq m-1}$ on $F$ such that 
\begin{equation}
\label{eq4}
\max_{|\alpha|\leq m-1}\Vert u^\alpha\Vert_{\infty,F}+
\max_{|\alpha|\leq m-1}\sup_{\substack{x,y\in F\\ x\neq y}}\frac{|u^\alpha(y)-T_x^{m-1-|\alpha|}u^\alpha(y)|}{|x-y|^{m-|\alpha|}}\leq M.
\end{equation}
Then there is a function $U\in W^{m,\infty}(\bbbr^n)\subset C^{m-1,1}(\bbbr^n)$ such that
$D^\alpha U=u^\alpha$ on $F$ and 
\begin{equation}
\label{eq12}
n^{-1}\Vert U\Vert_{m,\infty}\leq 
\sum_{|\alpha|\leq m-1}\Bigg(\Vert D^\alpha U\Vert_{\infty}+
\sup_{\substack{x,y\in\bbbr^n\\ x\neq y}}\frac{|D^\alpha U(y)-T_x^{m-1-|\alpha|}D^\alpha U(y)|}{|x-y|^{m-|\alpha|}}\Bigg)\leq C_wM,
\end{equation}
where the constant $C_w$ depends on $n$ and $m$ only, and the norms in \eqref{eq12} are over $\bbbr^n$.
\end{lemma}
\begin{remark}
The proof of Lemma~\ref{T3} presented in \cite{stein} is with $\max_{|\alpha|\leq m-1}$ in \eqref{eq12} instead of $\sum_{|\alpha|\leq m-1}$ (just like we have it in \eqref{eq4}), but the two quantities are clearly equivalent. The reason why we used $\sum_{|\alpha|\leq m-1}$ was to have a nice explicit constant in the left inequality in~\eqref{eq12}.
\end{remark}
\begin{remark}
If $m=1$, we simply have $u^0$ and \eqref{eq4} means that $u^0$ is bounded and Lipschitz continuous
$$
\sup_{x\in F} |u^0(x)|+\sup_{\substack{x,y\in F\\ x\neq y}}\frac{|u^0(y)-u^0(x)|}{|x-y|}\leq M.
$$
In that case the theorem says that $u^0$ can be extended to a bounded Lipschitz continuous function on $\bbbr^n$. 
\end{remark}
\begin{remark}
The proof given in \cite{stein}
shows the right inequality in \eqref{eq12}, but as we will see, the left inequality is nearly obvious.
Since $\sum_{|\alpha|\leq m-1}\Vert D^\alpha U\Vert_{\infty}$ is a part of the middle expression in \eqref{eq12}, 
it suffices to concentrate on the $L^\infty$ norm of derivatives of order $m$. For $\alpha$ satisfying $|\alpha|=m-1$ we have
\[
\begin{split}
n^{-1}\sum_{i=1}^n \Big\Vert\frac{\partial}{\partial x_i}D^\alpha U\Big\Vert_{\infty}
&\leq
\Vert\nabla D^\alpha U\Vert_{\infty}=
\sup_{\substack{x,y\in\bbbr^n\\ x\neq y}}\frac{|D^\alpha U(y)-D^\alpha U(x)|}{|x-y|} \\
&=
\sup_{\substack{x,y\in\bbbr^n\\ x\neq y}}\frac{|D^\alpha U(y)-T_x^{m-1-|\alpha|}D^\alpha U(y)|}{|x-y|^{m-|\alpha|}}\, ,
\end{split}
\]
and the left inequality in \eqref{eq12} follows.
\end{remark}
It should not be surprising now that Corollary~\ref{T2} and Lemma~\ref{T3} can be used to prove that $W^{m,1}$ functions coincide with $C^m$ functions outside a set of small measure.
This is a theorem of
Calder\'on and Zygmund \cite{CZ} who proved it using different techniques. Several generalizations of their result were obtained by many authors including \cite{BH,BHS,Liu,MichZie}. It seems the papers \cite{BH,BHS} were the first to use pointwise inequalities as in Lemma~\ref{T1} and Corollary~\ref{T2} in that context. A modified approach was given by Mal\'y and Ziemer \cite[Theorem~1.69]{MZ} who proved it for $m=1$ using pointwise inequalities from \cite{BH} (see \cite[p. 62]{MZ}). The next result is a generalization of  \cite[Theorem~1.69]{MZ} to the case of higher order derivatives. Recall that $W^m=W^{m,1}\cap W^{m,\infty}$.
\begin{proposition}
\label{T4}
Let $u\in W^{m,1}(\bbbr^n)$. Then, for every $a>0$, there is $v\in W^m\cap C^m(\bbbr^n)$, such that
\begin{align}
&\Vert v\Vert_{m,\infty}\leq a \label{eq5}\\
&\Vert v\Vert_{m,1}\leq C_*\Vert u\Vert_{m,1} \label{eq6}\\
&|\{x\in\bbbr^n:\, u(x)\neq v(x)\}|\leq\frac{C_*}{a}\Vert u\Vert_{m,1}, \label{eq7}
\end{align}
where the constant $C_*$ depends on $n$ and $m$ only.
\end{proposition}
\begin{proof}
Let $u\in W^{m,1}(\bbbr^n)$. In the first step, we will prove a weaker result, existence of a function $w\in W^m(\bbbr^n)$ that has similar properties
\begin{align}
&\Vert w\Vert_{m,\infty}\leq a \label{eq8}\\
&\Vert w\Vert_{m,1}\leq C_o\Vert u\Vert_{m,1} \label{eq9}\\
&|\{x\in\bbbr^n:\, u(x)\neq w(x)\}|\leq\frac{C_o}{a}\Vert u\Vert_{m,1}, \label{eq10}
\end{align}
where the constant $C_o$ depends on $n$ and $m$ only.
That is, we will construct a function with all properties required by the proposition except for the $C^m$ regularity.

Without loss of generality, we may assume that $u$ and its all derivatives are defined at every point of $\bbbr^n$ by \eqref{eq1}.

In view of Corollary~\ref{T2}, in order to apply Lemma~\ref{T3} with $u^\alpha:=D^\alpha u$, we need to restrict $u$ to some closed set $F$, where
$\max_{|\alpha|\leq m-1}|D^\alpha u|$ and $\MM|\nabla^m u|$ are bounded by suitable constants.
The set $F$ will be defined as the complement of the open set $G_0\cup G_1$ defined below.

Recall that the constants $C_p$ and $C_w$ from Corollary~\ref{T2} and Lemma~\ref{T3} depend on $n$ and $m$ only. We define
$$
\widetilde{G}_0:=\Big\{x\in\bbbr^n:\, \max_{|\alpha|\leq m-1}|D^\alpha u(x)|>\frac{a}{2nC_w}\Big\},
$$
$$
G_1:=\Big\{x\in\bbbr^n:\, \MM |\nabla^m u|(x)>\frac{a}{4nC_pC_w}\Big\}.
$$
Chebyshev's inequality and the weak type estimate for the maximal function \eqref{eq11} yield
$$
|\widetilde{G}_0|\leq\frac{2nC_w}{a}\, \Vert u\Vert_{m,1}
\quad
\text{and}
\quad
|G_1|\leq 5^n\frac{4nC_pC_w}{a}\Vert u\Vert_{m,1},
$$
Note that the set $G_1$ is open, but $\widetilde{G}_0$ not necessarily. Let $G_0$ be an open set such that $\widetilde{G}_0\subset G_0$ and
$$
|G_0|\leq\frac{4nC_w}{a}\Vert u\Vert_{m,1}.
$$
Let $G:=G_0\cup G_1$ and $F=G^c$. Then
\begin{equation}
\label{eq13}
|G|\leq 4nC_w(1+5^nC_p)\,\frac{1}{a}\,\Vert u\Vert_{m,1}.
\end{equation}
If $x\in F$, then we have inequalities opposite to those in the definition of the sets $\widetilde{G}_0$ and $G_1$, which combined with \eqref{eq3} yield
\[
\begin{split}
\max_{|\alpha|\leq m-1}\sup_{x\in F} |D^\alpha u(x)|
&+
\max_{|\alpha|\leq m-1}\sup_{\substack{x,y\in F\\ x\neq y}}\frac{|D^\alpha u(y)-T_x^{m-1-|\alpha|}D^\alpha u(y)|}{|x-y|^{m-|\alpha|}} \\
&\leq 
\frac{a}{2nC_w}+2C_p\, \frac{a}{4nC_pC_w}=\frac{a}{nC_w}.
\end{split}
\]
It follows from this estimate that the functions $u^\alpha:=D^\alpha u|_F$ are continuous.

Applying Lemma~\ref{T3} we get $w:=U\in W^{m,\infty}(\bbbr^n)$ satisfying $\Vert w\Vert_{m,\infty}\leq a$ which is \eqref{eq8}. Since $u=w$ in the complement of $G$, inequality \eqref{eq13} yields
$$
|\{x\in\bbbr^n:\, u(x)\neq w(x)\}|
\leq |G| \leq \frac{C(n,m)}{a}\Vert u\Vert_{m,1}
$$
which is \eqref{eq10}. Finally \eqref{eq9} follows from the estimate
$$
\Vert w\Vert_{m,1}
=
\sum_{|\alpha|\leq m}\left(\Vert D^\alpha w\Vert_{1,G}+\Vert D^\alpha w\Vert_{1,G^c}\right)
\leq 
\Vert w\Vert_{m,\infty} |G| +\Vert u\Vert_{m,1} \\
\leq 
C(n,m)\Vert u\Vert_{m,1}.
$$
In the second to last inequality we used the fact that $w=u$ in $G^c$, and in
the last inequality we used \eqref{eq8} and \eqref{eq13}.
This completes the construction of a function $w$ with properties \eqref{eq8}, \eqref{eq9} and \eqref{eq10}.

Now we are ready to complete the proof of the proposition.

Using the construction above (with $a$ replaced by $2^{-1}a$), we can find $w_1\in W^m(\bbbr^n)$ such that
\begin{equation}
\label{eq16}
\Vert w_1\Vert_{m,\infty}\leq2^{-1}a,
\qquad
\Vert w_1\Vert_{m,1}\leq C_o\Vert u\Vert_{m,1},
\end{equation}
\begin{equation}
\label{eq17}
|\underbrace{\{x\in\bbbr^n:\, u(x)\neq w_1(x)\}}_{:=E_1}|\leq 2C_o\, a^{-1}\Vert u\Vert_{m,1}.
\end{equation}
Let $\vi$ be a mollifier as in Section~\ref{S1}. Then according to \eqref{eq14} and \eqref{eq15}, there is $\eps>0$ such that
\begin{equation}
\label{eq21}
v_1:=\vi_\eps*w_1
\qquad
\text{and}
\qquad
u_1:=w_1-v_1
\end{equation}
satisfy
\begin{align}
&\Vert v_1\Vert_{m,\infty}=\Vert\vi_\eps*w_1\Vert_{m,\infty}\leq\Vert w_1\Vert_{m,\infty}\leq2^{-1}a,\label{eq19}\\
&\Vert v_1\Vert_{m,1}=\Vert\vi_\eps*w_1\Vert_{m,1}\leq\Vert w_1\Vert_{m,1}\leq C_o\Vert u\Vert_{m,1},\label{eq20}\\
&\Vert u_1\Vert_{m,1}=\Vert w_1-\vi_\eps*w_1\Vert_{m,1}\leq 4^{-1}\Vert u\Vert_{m,1},\label{eq18}\\
&\Vert u_1\Vert_\infty\leq\Vert w_1\Vert_\infty+\Vert v_1\Vert_\infty\leq a. \label{eq22}
\end{align}

The function $v_1$ is $C^\infty$ smooth and it satisfies \eqref{eq19} and \eqref{eq20} which yield $v_1\in W^m\cap C^\infty(\bbbr^m)$, \eqref{eq5} and \eqref{eq6}, but the problem is that $v_1$ need not coincide with $u$ on a set of positive measure. 
On the other hand
$$
v_1+u_1=u
\quad 
\text{in $E_1^c$}
$$
by \eqref{eq21} and \eqref{eq17}, but the `correcting term' $u_1$ is not smooth. It belongs to $W^{m,1}$, but it has much smaller $W^{m,1}$ norm than $u$, see \eqref{eq18}.
The idea is to apply the above procedure to $u_1$ in place of $u$ and find $v_2$, $u_2$, and $E_2$, so that $v_2+u_2=u_1$ in $E_2^c$ and hence
$$
v_1+v_2+u_2=v_1+u_1=u 
\quad
\text{in $E_1^c\cap E_2^c$.}
$$
Continuing this procedure, we will have
$$
u_k+\sum_{i=1}^k v_i=u
\quad
\text{in}
\quad
\bigcap_{i=1}^k E_i^c.
$$
It will follow from the estimates that the series $\sum_{i=1}^\infty v_i$ will converge to a function $v\in W^m\cap C^m(\bbbr^n)$, and $u_k$ will converge uniformly to zero, so $v=u$ in $\bigcap_{i=1}^\infty E_i^c$. That will complete the proof.

The sequence of functions and sets will be constructed by induction. 
We already constructed the functions $v_1\in W^m\cap C^\infty$, $u_1\in W^{m,1}$ and the set $E_1$. Suppose now that for some $k\in\bbbn$ we have functions $\{v_i\}_{i=1}^k\subset W^m\cap C^\infty(\bbbr^n)$, $\{u_i\}_{i=1}^k\subset W^{m,1}(\bbbr^n)$, and sets $\{E_i\}_{i=1}^k$ such that 
for $i=1,2,\ldots,k$ we have
\begin{align}
&\Vert v_i\Vert_{m,\infty}\leq 2^{-i}a, \label{eq23}\\
&\Vert v_i\Vert_{m,1}\leq 4^{-i+1}\, C_o\Vert u\Vert_{m,1}, \label{eq24}\\
&\Vert u_k\Vert_{m,1}\leq 4^{-k}\Vert u\Vert_{m,1}, \label{eq25}\\
&\Vert u_k\Vert_\infty\leq 2^{-k+1}a, \label{eq26}\\
&u_k+\sum_{i=1}^k v_i=u
\quad
\text{in} 
\quad
\bigcap_{i=1}^k E_i^c, \label{eq27}\\
&|E_i|\leq 2^{-i+2}\, C_oa^{-1}\Vert u\Vert_{m,1}. \label{eq28}
\end{align}
We already verified these conditions for $k=1$.

Now we will construct $v_{k+1}$, $u_{k+1}$ and $E_{k+1}$ and we will verify conditions \eqref{eq23}--\eqref{eq28}.

Applying the construction of $w$ satisfying \eqref{eq8}, \eqref{eq9}, \eqref{eq10} to $u_k$ and $2^{-(k+1)}a$ (instead of $u$ and $a$), we find $w_{k+1}\in W^m(\bbbr^n)$ such that
$$
\Vert w_{k+1}\Vert_{m,\infty}\leq 2^{-(k+1)}a,
\qquad
\Vert w_{k+1}\Vert_{m,1}\leq C_o\Vert u_k\Vert_{m,1}\leq 4^{-k}\, C_o\Vert u\Vert_{m,1},
$$
$$
|\underbrace{\{x\in\bbbr^n:\, u_k(x)\neq w_{k+1}(x)\}}_{:=E_{k+1}}|\leq
2^{k+1}\, C_o\,a^{-1}\Vert u_k\Vert_{m,1}\leq
2^{-(k+1)+2}\, C_oa^{-1}\Vert u\Vert_{m,1}
$$
which proves \eqref{eq28} for $i=k+1$.

Now there is $\eps>0$ (perhaps different than the one before), such that
$$
v_{k+1}=\vi_\eps*w_{k+1}
\qquad
\text{and}
\qquad
u_{k+1}=w_{k+1}-v_{k+1}
$$
satisfy
\begin{align*}
& \Vert v_{k+1}\Vert_{m,\infty}\leq \Vert w_{k+1}\Vert_{m,\infty}\leq 2^{-(k+1)}\, a,\\
& \Vert v_{k+1}\Vert_{m,1}\leq\Vert w_{k+1}\Vert_{m,1}\leq 4^{-(k+1)+1}\, C_o\Vert u\Vert_{m,1},\\
& \Vert u_{k+1}\Vert_{m,1}=\Vert w_{k+1}-\vi_\eps*w_{k+1}\Vert_{m,1}\leq 4^{-(k+1)}\Vert u\Vert_{m,1},\\
&\Vert u_{k+1}\Vert_\infty\leq \Vert w_{k+1}\Vert_\infty+\Vert v_{k+1}\Vert_\infty\leq 2^{-(k+1)+1}\, a.
\end{align*}
This proves \eqref{eq23}, \eqref{eq24}, \eqref{eq25} and \eqref{eq26} for $i=k+1$. It remains to prove \eqref{eq27} for $i=k+1$.

Since $u_{k+1}+v_{k+1}=w_{k+1}$ and $w_{k+1}=u_k$ in $E_{k+1}^c$, \eqref{eq27} yields
$$
u_{k+1}+\sum_{i=1}^{k+1} v_i=
w_{k+1}+\sum_{i=1}^k v_i=
u_k+\sum_{i=1}^k v_i=u
\quad
\text{in}
\quad
\bigcap_{i=1}^{k+1} E_i^c.
$$
which proves \eqref{eq27} for $k+1$.

The proof of \eqref{eq23}--\eqref{eq28} for all $k\in\bbbn$ is complete and we are ready to complete the proof of the proposition.

Consider the series $v:=\sum_{i=1}^\infty v_i$. It follows from \eqref{eq23} and \eqref{eq24} that the series converges both in $W^{m,\infty}$ and $W^{m,1}$, so $v\in W^m(\bbbr^n)$ satisfies 
\begin{equation}
\label{eq29}
\Vert v\Vert_{m,\infty}\leq a
\qquad
\text{and}
\qquad
\Vert v\Vert_{m,1}\leq \frac{4}{3}\, C_o\Vert u\Vert_{m,1}.
\end{equation}
Since this is a series of smooth functions and $W^{m,\infty}\cap C^m(\bbbr^n)$ is a closed subspace of $W^{m,\infty}(\bbbr^n)$, it follows that
\begin{equation}
\label{eq30}
v\in W^m\cap C^m(\bbbr^n).
\end{equation}
Finally, \eqref{eq26} and \eqref{eq27} show that the series converges uniformly to $u$ on $\bigcap_{i=1}^\infty E_i^c$, so $v=u$ in that set and hence \eqref{eq28} yields
\begin{equation}
\label{eq31}
|\{x\in\bbbr^n:\, u(x)\neq v(x)\}|\leq\Big|\bigcup_{i=1}^\infty E_i\Big|\leq 4C_oa^{-1}\Vert u\Vert_{m,1}.
\end{equation}
Now the proposition follows from \eqref{eq29}, \eqref{eq30} and \eqref{eq31} with $C_*=4C_o$. The proof is complete.
\end{proof}

\section{Proof of Theorem~\ref{T22}}
\label{S3}

First, we will prove a special case of the result.
\begin{lemma}
\label{T25}
Let $u:\bbbr^n\to\bbbr$ be $\eta$-strongly convex. Then for every $\eps>0$ and every $0<\widetilde{\eta}<\eta$, there is an $\widetilde{\eta}$-strongly convex function 
$v\in C^2(\bbbr^n)$, such that
$$
|\{x\in\bbbr^n:\, u(x)\neq v(x)\}|<\eps.
$$
\end{lemma}
\begin{proof}
We claim it is enough to prove the result for $\eta$-strongly convex functions in $W^{2,p}_{\rm loc}(\bbbr^n)$ for $n<p<\infty$. Indeed, the general case can be concluded then as follows:

Suppose the result holds for functions in $W^{2,p}_{\rm loc}(\bbbr^n)$ for $n<p<\infty$. 
Let $u:\bbbr^n\to\bbbr$ be $\eta$-strongly convex, and fix $\eps>0$ and $0<\widetilde{\eta}<\eta$.
Let $\eta' \in (\widetilde{\eta},\eta)$. We apply Corollary~\ref{T11} to produce ${u'}\in C^{1,1}_{\rm loc}$ that it $\eta'$-strongly convex and satisfies 
$|\{ {u'}\neq u\}|<\eps/2$. In particular, ${u'}\in W^{2,p}_{\rm loc}$ for any $n<p<\infty$. 
Then, we apply the assumed result to ${u'}$, $\eps/2>0$, and $\widetilde{\eta} \in (0,\eta')$, producing $v\in C^2$ that is $\widetilde{\eta}$-strongly convex and satisfies $|\{ {u'}\neq v\}|<\eps/2$, and thus $|\{u\neq v\}|<\eps$.

Therefore, in what follows we assume  that $u\in W^{2,p}_{\rm loc}(\bbbr^n)$ for some $n<p<\infty$, is $\eta$-strongly convex. Without loss of generality, we may assume that
$$
0<\eps<\frac{\eta-\widetilde{\eta}}{2nC_*},
$$
where $C_*$ is the constant from Proposition~\ref{T4}.

According to Proposition~\ref{T6}, there is $\widetilde{u}\in C^\infty(\bbbr^n)$ such that $\widetilde{u}$ is $(\eta+\widetilde{\eta})/2$-strongly convex and $\Vert u-\widetilde{u}\Vert_{2,1}<\eps^2$.

Applying Proposition~\ref{T4} with $m=2$ and $a=C_*\eps$ to $u-\widetilde{u}$ in place of $u$, we find $w\in W^2\cap C^2(\bbbr^n)$ such that
\begin{equation}
\label{eq36}
\Vert w\Vert_{2,\infty}\leq C_*\eps<\frac{\eta-\widetilde{\eta}}{2n},
\end{equation}
$$
\Vert w\Vert_{2,1}\leq C_*\Vert u-\widetilde{u}\Vert_{2,1}<C_*\eps^2,
$$
$$
\Big|\big\{x\in\bbbr^n:\, w(x)\neq u(x)-\widetilde{u}(x)\big\}\big|\leq\frac{C_*}{C_*\eps}\, \Vert u-\widetilde{u}\Vert_{2,1}<\eps.
$$
Since 
$$
\widetilde{u} 
\text{ is } \frac{\eta+\widetilde{\eta}}{2}\text{--strongly convex}, 
\qquad
0<\widetilde{\eta}<\frac{\eta+\widetilde{\eta}}{2},
$$
and according to \eqref{eq36}
$$
n\max_{|\alpha|=2}\Vert D^\alpha w\Vert_\infty\leq\frac{\eta-\widetilde{\eta}}{2}=\frac{\eta+\widetilde{\eta}}{2}-\widetilde{\eta},
$$
it follows from Lemma~\ref{T12} that $v:=\widetilde{u}+w$ is $\widetilde{\eta}$-strongly convex. Clearly, $v\in C^2(\bbbr^n)$ and
$|\{v\neq u\}|=|\{w\neq u-\widetilde{u}\}|<\eps$.
The proof is complete.
\end{proof}

\begin{corollary}
\label{T19}
Let $u:\R^n\to\R$ be $\eta$-strongly convex, $\widetilde{\eta}\in (0, \eta)$, and let $0<\eps_k\to 0$. Then, there exists a sequence of $\widetilde{\eta}$-strongly convex functions $(v_k)_{k \in \N} \subset C^2(\R^n)$ satisfying:
\begin{enumerate}
\item[(a)] $|\{x\in\R^n: v_k(x)\neq u(x)\}|<\varepsilon_k$ for all $k\in\N$;
\item[(b)] $v_k\rightrightarrows u$ converges uniformly on every compact subset of $\R^n$.
\end{enumerate}
\end{corollary}
\begin{proof}
Given $u, \eta, \widetilde{\eta}$, and $(\varepsilon_k)_{k \in \N}$ satisfying the hypotheses of the corollary, we apply Lemma~\ref{T25} to obtain a sequence  of $\widetilde{\eta}$-strongly convex functions $(v_k)_{k \in \N} \subset C^2(\R^n)$, satisfying $|\{x\in\R^n: v_k(x)\neq u(x)\}|<\varepsilon_k/2^k$ for all $k\in\N$ (and consequently (a)). 

It suffices to show that $v_k\to u$ almost everywhere, because then (b) will follow from Lemma~\ref{T14}.
Let $A_i:=\{x\in\R^n:\, v_i(x)=u(x)\}$, and $C:=\bigcup_{k=1}^\infty\bigcap_{i=k}^\infty A_i$. If $x\in C$, then $x\in\bigcap_{i=k}^\infty A_i$ for some $k$ and hence
$u(x)=v_k(x)=v_{k+1}(x)=\ldots$,
so
$v_i(x)\stackrel{i\to\infty}{\longrightarrow}u(x)$.
We proved that $v_i(x)\to u(x)$ for every $x\in C$ and it remains to show that $|\R^n\setminus C|=0$. Clearly,
$|\R^n\setminus A_i|<\eps_i/2^i$. Since
$$
|\R^n\setminus C|\leq \Big| \R^n\setminus\bigcap_{i=k}^\infty A_i\Big|=
\Big|\bigcup_{i=k}^\infty (\R^n\setminus A_i)\Big|
<\sum_{i=k}^\infty\frac{\eps_i}{2^i}\leq \sup_{i\geq k}\eps_i\stackrel{k\to\infty}{\longrightarrow} 0,
$$
we conclude that $|\R^n\setminus C|=0$.
\end{proof}

\begin{proof}[Proof of Theorem~\ref{T22}]
Let $u:U\to\R$ be locally strongly convex, $\eps_o>0$, and $\varepsilon:U\to (0,1]$ continuous. 
We fix a sequence $(B_k)_{k=1}^\infty$ of compact convex bodies such that
$B_k\subset\inte(B_{k+1})$ for all $k\in\bbbn$ and $\bigcup_{k=1}^\infty B_k=U$.

By assumption, for each $k\in\bbbn$, there is $\eta_k>0$ such that $u|_{B_k}$ is $\eta_k$-strongly convex; we may assume that $\eta_{k+1}<\eta_k$ for every $k$.

Let
\begin{equation}
\label{eq52}
r_k:=\min_{x\in B_{k+1}}\eps(x)
\quad
\text{for } k\in\N.
\end{equation}
Let $B_0=C_0=\varnothing$. For each $k\geq 1$, we find a compact convex body $C_k$ such that
\begin{equation}
\label{eq48}    
\text{$B_k\subset \inte(C_k)\subset C_k\subset \inte(B_{k+1})$}, 
\quad
\dist(\partial B_k,\partial C_k)\leq\frac{r_{k+1}}{6\Lip(u|_{B_{k+1}})},
\end{equation}
and
\begin{equation}
\label{eq44}
|C_k \setminus B_k|\leq \frac{\eps_o}{2^{k+2}\sqrt{1+L_{k+1}^2}}, 
\quad
\text{where}
\quad
L_k:=\frac{2+M_{B_{k+1}}(u)-m_{B_{k+1}}(u)}{\dist(\partial B_k, \partial B_{k+1})}
\end{equation}
(we are using notation from Lemma~\ref{T13}). Note that $\Lip(u|_{B_{k+1}})>0$, because $u|_{B_{k+1}}$ is strongly convex, so we do not divide by zero in \eqref{eq48}.

We can also assume that
\begin{equation}
\label{eq53}
C_k=\{x:\dist(x,B_k)\leq a_k\}
\end{equation}
for some $a_k>0$ so that all points on the boundary of $C_k$ are at the equal distance to $B_k$. Clearly, $a_k=\dist(\partial B_k,\partial C_k)$.

Note that \eqref{eq44} implies
\begin{equation}
\label{eq45}
\sum_{k=1}^{\infty}|C_k\setminus B_k|<\frac{\eps_o}{2}.
\end{equation}
Then, for each $k\in\N$ we define $u_k:\R^n\to\R$ by (cf.\ Corollary~\ref{T18})
$$
u_{k}(x):=\sup_{y\in B_{k}\setminus \inte(C_{k-1}), \, \xi_y\in\partial u(y)}\left\{u(y)+\langle \xi_y, x-y\rangle +\frac{\eta_{k+8}}{2}|x-y|^2\right\}.
$$
The next result collects important properties of the functions $u_k$.
\begin{lemma}
\label{T20}
For each $k\in\N$, the function $u_k:\R^n\to\R$ satisfies:
\begin{enumerate}
\item[(a)] $u_k$ is $\eta_{k+8}$-strongly convex;
\item[(b)] $u_k=u$ on $B_k\setminus C_{k-1}$;
\item[(c)] $u_k\leq u$ on $B_{k+8}$;
\item[(d)] $\beta_k:=\inf_{x\in B_k}\{u(x)-u_{k+1}(x)\}>0$; 
\item[(e)] $\delta_k:=\inf_{x\in B_{k+1}\setminus \inte(C_{k})}\{u(x)-u_k(x)\}>0$;
\item[(f)] $\Vert u_k-u\Vert_{\infty,C_k\setminus B_{k-1}}\leq r_k/3$; recall that $B_0=\varnothing$.
\end{enumerate}
\end{lemma}
\begin{proof}
(a), (b) and (c) follow from Corollary~\ref{T18}. To prove (d) let $x\in B_k$ and $y\in B_{k+1}\setminus\inte(C_k)$, $\xi_y\in \partial u(y)$. Since $u$ is $\eta_{k+2}$-strongly convex on $B_{k+2}$, 
and $x,y\in\inte(B_{k+2})$,
Lemma~\ref{T16} yields
\begin{align*}
u(x)&\geq u(y) +\langle \xi_y, x-y\rangle +\frac{\eta_{k+2}}{2}|x-y|^2 \\
&\geq u(y) +\langle \xi_y, x-y\rangle +\frac{\eta_{k+9}}{2}|x-y|^2 +\frac{\eta_{k+2}-\eta_{k+9}}{2}{\dist(\partial B_k, \partial C_k)}^2.
\end{align*}
Taking the supremum over all $y\in B_{k+1}\setminus\inte(C_k)$ and $\xi_y\in \partial u(y)$ gives
$$
u(x)\geq u_{k+1}(x)+\frac{\eta_{k+2}-\eta_{k+9}}{2}{\dist(\partial B_k, \partial C_k)}^2
\quad
\text{for all } x\in B_k,
$$
from which (d) follows. The proof of (e) is similar to that of (d). For $x\in B_{k+1}\setminus\inte(C_k)$, $y\in B_k\setminus\inte(C_{k-1})$ and $\xi_y\in\partial u(y)$, $\eta_{k+2}$-strong convexity of $u$ on $B_{k+2}$ and Lemma~\ref{T16} yield
\begin{align*}
u(x)&\geq u(y) +\langle \xi_y, x-y\rangle +\frac{\eta_{k+2}}{2}|x-y|^2 \\
&\geq u(y) +\langle \xi_y, x-y\rangle +\frac{\eta_{k+8}}{2}|x-y|^2 +\frac{\eta_{k+2}-\eta_{k+8}}{2}{\dist(\partial B_k, \partial C_k)}^2.
\end{align*}
Taking the supremum over all $y\in B_k\setminus\inte(C_{k-1})$ and $\xi_y\in\partial u(y)$ gives
$$
u(x)\geq u_{k}(x)+\frac{\eta_{k+2}-\eta_{k+8}}{2}{\dist(\partial B_k, \partial C_k)}^2
\quad
\text{for all } x\in B_{k+1}\setminus\inte(C_k),
$$
from which (e) follows. 

It remains to prove (f). Since $u=u_k$ on $B_k\setminus C_{k-1}$, we only need to consider the sets $C_{k-1}\setminus B_{k-1}$ and $C_k\setminus B_k$, when $k\geq 2$, and $C_k\setminus B_k$, when $k=1$.

If $x\in C_{k-1}\setminus B_{k-1}$, $k\geq 2$, then we can find $y\in\partial C_{k-1}\subset B_k\setminus\inte(C_{k-1})$ such that
\begin{equation}
\label{eq54}
|x-y|\leq\dist(y,\partial B_{k-1})=\dist(\partial B_{k-1},\partial C_{k-1})
\end{equation}
(in the last equality we are using \eqref{eq53}). If $\xi_y\in\partial u(y)$, then
$|\xi_y|\leq\Lip(u|_{B_k})$ by Lemma~\ref{T23}, because $y\in\inte(B_k)$. Property (c) of Lemma~\ref{T20} yields
\[
\begin{split}
u(x)
&\geq
u_k(x)\geq u(y)+\langle\xi_y,x-y\rangle +\frac{\eta_{k+8}}{2}|x-y|^2\\
&\geq
u(x)-|u(y)-u(x)|-|\xi_y|\, |x-y|\\
&\geq 
u(x)-2\Lip(u|_{B_k})|x-y|\geq u(x)-\frac{r_k}{3}.
\end{split}
\]
In the last inequality we used \eqref{eq54} and \eqref{eq48}. Therefore,
$$
\Vert u_k-u\Vert_{\infty,C_{k-1}\setminus B_{k-1}}\leq\frac{r_k}{3}
\quad
\text{for $k\geq 2$.}
$$
If $x\in C_k\setminus B_k$, $k\geq 1$, then we can find $y\in\partial B_k\subset B_k\setminus\inte(C_{k-1})$ such that
$|x-y|\leq\dist(\partial B_k,\partial C_k)$. If $\xi_y\in\partial u(y)$, then $|\xi_y|\leq \Lip(u|_{B_{k+1}})$ and 
\[
\begin{split}
u(x)
&\geq 
u_k(x)\geq u(x)-|u(y)-u(x)|-|\xi_y|\, |x-y|\\
&\geq 
u(x)-2\Lip(u|_{B_{k+1}})|x-y|\geq u(x)-\frac{r_{k+1}}{3}.
\end{split}
\]
Therefore,
$$
\Vert u_k-u\Vert_{\infty,C_k\setminus B_k}\leq\frac{r_{k+1}}{3}\leq\frac{r_k}{3}\, .
$$
This completes the proof of (f).
\end{proof}
The idea of the remaining part of the proof is to use Corollary~\ref{T19} to approximate $u_k$ near the annulus  $B_k\setminus\inte(C_{k-1})$ by globally defined strongly convex functions of class $C^2$ and glue these approximations using the smooth maximum method described in Lemma~\ref{T17}.

Let $\eps_k$ be a sequence such that
$$
0<\eps_k\leq \frac{1}{3}
\min\left\{\frac{\min\{\eps_o,r_k\}}{2^{k+2} \sqrt{1+ L_k^2}}, \delta_k, \beta_k\right\},
\quad
\eps_{k+1}\leq\eps_k
\quad
\text{for } k\in\N
$$
(recall that $r_k$ and $L_k$ were defined in \eqref{eq52} and \eqref{eq44}).
In particular,
\begin{equation}
\label{eq50}
\sum_{k=1}^\infty \eps_k<\frac{\eps_o}{2}.
\end{equation}

Now, we use
Corollary~\ref{T19} to find an $\eta_{k+9}$-strongly convex function $h_k\in C^{2}(\R^n)$ such that
\begin{equation}
\label{eq58}
|\{x\in\R^n : u_k(x)\neq h_k(x)\}|<\varepsilon_k
\end{equation}
and
\begin{equation}
\label{eq46}
|u_k(x)-h_k(x)|\leq\varepsilon_k \textrm{ for all } x\in B_{k+1}.
\end{equation}
For $k\geq 2$, we define  $v_{k}:\R^n\to\R$ by (see, \eqref{eq57})
$$
v_k=M_{\varepsilon_{k-1}}(h_{k-1}, h_{k}).
$$
\begin{lemma}
\label{T21}
For each $k\geq 2$, the function $v_k$ is in the class $C^2(\R^n)$, $\eta_{k+9}$-strongly convex, and satisfies:
\begin{enumerate}
\item[(a)] $v_k=h_{k-1}$ on $B_{k-1}\setminus C_{k-2}$;
\item[(b)] $v_{k}=h_{k}$ on $B_{k}\setminus C_{k-1}$.
\end{enumerate}
\end{lemma}
\begin{proof}
The fact that $v_k\in C^2(\R^n)$ and $v_k$ is $\eta_{k+9}$-strongly convex, follows immediately from Lemma~\ref{T17}.

\noindent 
$(a)$ 
For $x\in B_{k-1}\setminus C_{k-2}$, \eqref{eq46} and Lemma~\ref{T20} give
\begin{align*}
h_{k}&\leq u_{k}+\varepsilon_{k}\leq u-\beta_{k-1}+\varepsilon_{k}
=u_{k-1}-\beta_{k-1}+\varepsilon_{k}\\
&\leq h_{k-1} +\varepsilon_{k-1}-\beta_{k-1} +\varepsilon_{k}
\leq h_{k-1}-\beta_{k-1} +2\varepsilon_{k-1},
\end{align*}
hence 
$$
h_{k-1}-h_{k}\geq \beta_{k-1}-2\varepsilon_{k-1}\geq \varepsilon_{k-1},
$$
and from Lemma~\ref{T17}(b) we deduce that $v_{k}=M_{\varepsilon_{k-1}}(h_{k-1}, h_{k})=h_{k-1}$ on $B_{k-1}\setminus C_{k-2}$.

\noindent
$(b)$
For $x\in B_{k}\setminus C_{k-1}$,
\eqref{eq46} and Lemma~\ref{T20} give
\begin{align*}
h_{k-1}&\leq u_{k-1} +\varepsilon_{k-1}\leq u-\delta_{k-1} +\varepsilon_{k-1}=u_{k}-\delta_{k-1}+\varepsilon_{k-1}\\
&\leq h_{k} +\varepsilon_{k}-\delta_{k-1}+\varepsilon_{k-1}\leq h_{k}-\delta_{k-1}+2\varepsilon_{k-1},
\end{align*}
hence
$$
h_{k}-h_{k-1}\geq\delta_{k-1}-2\varepsilon_{k-1}\geq\varepsilon_{k-1},
$$
and we conclude that $v_k=M_{\varepsilon_{k-1}}(h_{k-1}, h_{k})=h_{k}$ on $B_{k}\setminus C_{k-1}$.
\end{proof}
Now, we define our function $v:U\to\R$ by
\begin{equation}
\label{eq47}
v(x)=v_{k}(x) \text{ if } x\in B_{k}\setminus C_{k-2} \text{ for some } k\geq 2.
\end{equation}
While $\bigcup_{k=2}^\infty(B_k\setminus C_{k-2})=U$, consecutive annuli in the definition \eqref{eq47} overlap. However, for $k\geq 2$, on the overlapping annuli we have
\begin{equation}
\label{eq51}
v=v_{k}=v_{k+1}=h_k
\quad
\text{on } 
\quad
(B_k\setminus C_{k-2}) \cap (B_{k+1}\setminus C_{k-1})=B_k\setminus C_{k-1},
\end{equation}
showing that $v$ is well defined on $U$. 
Moreover, since the functions $v_k\in C^2(\R^n)$ are strongly convex  and $v_{k+1}=v_k$ on $B_k\setminus C_{k-1}$ for every $k\geq 2$, the function $v$ is locally strongly convex and of class $C^2(U)$. 

Note that we have
\begin{equation}
\label{eq59}
v=h_k 
\quad
\text{on}
\quad
B_k\setminus C_{k-1}
\quad
\text{for all } k\in\N.
\end{equation}
If $k\geq 2$, it follows from \eqref{eq51}, and we check the case $k=1$ directly from \eqref{eq47} and Lemma~\ref{T21}(a).

\noindent
({\bf Remark.} If all functions $v_k$ are $\widetilde{\eta}$-strongly convex, then $v$ is $\widetilde{\eta}$-strongly convex; we will need this fact in the last part of the proof.)

Observe that
$$
U=\bigcup_{k=1}^\infty (C_k\setminus B_k)\, \cup\, \bigcup_{k=1}^\infty (B_k\setminus C_{k-1}).
$$
Moreover, $v=h_k$ in $B_k\setminus C_{k-1}$ by \eqref{eq59}, and $u=u_k$ in $B_k\setminus C_{k-1}$ by Lemma~\ref{T20}(b). Therefore,
\[
\begin{split}
|\{x\in U: v(x)\neq u(x)\}|
&\leq
\sum_{k=1}^{\infty}|C_k\setminus B_k|+
\sum_{k=1}^{\infty}|\{x\in B_{k}\setminus C_{k-1} :\, h_k(x)\neq u_k(x)\}| \\
&<
\frac{\eps_o}{2}+\sum_{k=1}^{\infty}\eps_k<\eps_o.
\end{split}
\]
In the last two inequalities we used \eqref{eq45}, \eqref{eq58} and \eqref{eq50}.

Thus, we constructed a locally strongly convex function $v\in C^2(U)$ satisfying condition (a) in Theorem~\ref{T22}.  Now we will prove (b) i.e.,we will prove that $|u(x)-v(x)|<\eps(x)$ for all $x\in U$.

If $x\in B_k\setminus C_{k-1}$ for some $k\in\N$, then Lemma~\ref{T20}(b), \eqref{eq59} and 
\eqref{eq46} give 
$$
|u(x)-v(x)|=|u_k(x)-h_k(x)|  \leq\eps_k <r_k\leq\eps(x).
$$

Thus, we may assume that $x\in C_k\setminus B_k$ for some $k\in\N$.
It follows from Lemma~\ref{T20}(f) that
$$
|u(x)-u_k(x)|\leq\frac{r_k}{3}
\quad
\text{and}
\quad
|u(x)-u_{k+1}(x)|\leq \frac{r_{k+1}}{3}\leq\frac{r_k}{3}.
$$
Also, \eqref{eq46} yields
$$
|u_k(x)-h_k(x)|\leq\eps_k<\frac{r_k}{3}
\quad
\text{and}
\quad
|u_{k+1}(x)-h_{k+1}(x)|\leq\eps_{k+1}<\frac{r_k}{3},
$$
so
\begin{equation}
\label{eq55}
|u(x)-h_k(x)|<\frac{2r_k}{3}
\quad
\text{and}
\quad
|u(x)-h_{k+1}(x)|<\frac{2r_k}{3}.
\end{equation}
Recall that $v=v_{k+1}=M_{\eps_k}(h_k,h_{k+1})$ on $C_k\setminus B_k$ and hence
$$
\max\{ h_k,h_{k+1}\}\leq v=M_{\eps_k}(h_k,h_{k+1})<\max\{ h_k,h_{k+1}\}+\frac{r_k}{3}
\quad
\text{on } C_k\setminus B_k
$$
by Lemma~\ref{T17}(c), because $\eps_k/2<r_k/3$. Therefore, \eqref{eq55} gives
$$
u(x)<h_k(x)+\frac{2r_k}{3}\leq v(x)+\frac{2r_k}{3}
$$
and
$$
u(x)>\max\{h_k(x),h_{k+1}(x)\}-\frac{2r_k}{3}>v(x)-r_k,
$$
so we have $|u(x)-v(x)|<r_k\leq\eps(x)$. This completes the proof of part (b) of the theorem.

Now we will show (c). According to Lemma~\ref{T13}, and the definition of $L_k$ in \eqref{eq44},
$$
\Lip(u|_{B_k})\leq\frac{M_{B_{k+1}}(u)- m_{B_{k+1}}(u)}{\dist(\partial B_k,\partial B_{k+1})}< L_k.
$$ 
In particular $|Du|\leq L_k$ almost everywhere in $B_k$. Recall also that
$$
u=u_k,\ 
v=h_k
\text{ in } B_k\setminus C_{k-1}
\quad
\text{and} 
\quad
|\{x\in\R^n:\, u_k(x)\neq h_k(x)\}|<\eps_k,
$$
see, Lemma~\ref{T20}(b), \eqref{eq59} and \eqref{eq58}.
Then the formula for the surface area of the graph yields
\[
\begin{split}
&\mathcal{H}^{n}\left(\mathcal{G}_{u}\setminus \mathcal{G}_{v}\right)=
\int_{\{x\in U: u(x)\neq v(x)\}}\sqrt{1+|Du(x)|^2}\,dx \\
&\leq 
\sum_{k=1}^\infty\int_{C_k\setminus B_k}\sqrt{1+L_{k+1}^2}\, dx +
\sum_{k=1}^\infty\int_{\{x\in B_k\setminus C_{k-1}:\, u_k(x)\neq h_k(x)\}}
\sqrt{1+L_{k}^2}\, dx\\
&\leq
\sum_{k=1}^{\infty}\sqrt{1+L_{k+1}^2}\, |C_k\setminus B_k| +\sum_{k=1}^{\infty} \sqrt{1+L_{k}^2}\,\varepsilon_k \\
& < \sum_{k=1}^{\infty}\frac{\eps_o}{2^{k+2}} +\sum_{k=1}^{\infty} \frac{\eps_o}{2^{k+2}}=\frac{\eps_o}{2}.
\end{split}
\]
On the other hand, since $|u(x)-v(x)|<\eps(x)\leq 1$ on $U$, we have
$$
\Lip(v|_{B_k})\leq \frac{M_{B_{k+1}}(v)-m_{B_{k+1}}(v)}{\textrm{dist}(\partial B_k, \partial B_{k+1})}\leq 
\frac{2+M_{B_{k+1}}(u)-m_{B_{k+1}}(u)}{\textrm{dist}(\partial B_k, \partial B_{k+1})}=L_k,
$$
and then the same calculation with $v$ in place of $u$ in the integrand shows that
$$
\mathcal{H}^{n}\left(\mathcal{G}_{v}\setminus \mathcal{G}_{u}\right)<\eps_o/2.
$$
Hence, we have $\mathcal{H}^{n}\left(\mathcal{G}_{u}\triangle \mathcal{G}_{v}\right)<\eps_o$.

Finally, in the case that $u$ is $\eta$-strongly convex on $U$, given $\widetilde{\eta}\in (0, \eta)$, we may find a strictly decreasing sequence $(\eta_k)_{k\in\N}\subset (\widetilde{\eta}, \eta)$ converging to $\widetilde{\eta}$ so that $u|_{B_k}$ is $\eta_k$-strongly convex, and repeat the proof above in order to obtain an $\widetilde{\eta}$-strongly convex function $v$ of class $C^2$ with the required properties (see Remark below \eqref{eq59}). 
\end{proof}

\end{document}